\documentclass[12pt]{amsart}
\usepackage{amsfonts}
\usepackage{amsmath,amssymb,amsthm}
\usepackage{amstext}
\usepackage{verbatim}
\input{epsf}
\newtheorem{thm}{Theorem}[section]

\newtheorem{prop}[thm]{Proposition}
\newtheorem{prob}[thm]{Problem}
\theoremstyle{definition}

\theoremstyle{remark}
\newtheorem{rem}[thm]{Remark}

\begin{document}

\title[Composition of inverse problems with a given logical structure]
{Composition of inverse problems with a given logical structure}
\author{Julia Ninova, Vesselka Mihova}%
\address{Faculty of Mathematics and Informatics
University of Sofia
e-mail: julianinova@hotmail.com}%
\address{Faculty of Mathematics and Informatics
University of Sofia
e-mail: mihova@fmi.uni-sofia.bg}


\maketitle
\thispagestyle{empty}
\vskip 4mm

\begin{abstract} The paper presents a method for obtaining
problems whose conclusions contain disjunctive propositions. These
problems constitute a version of inverse problems with a given
logical structure. The logical models in the groups of problems
studied have been interpreted comprehensively. Equivalent problems
have been given by keeping or not keeping the condition of
homogeneity in their conclusion. \vskip 4mm
\end{abstract}

\section{Introduction}

In mathematical logic a propositional calculus (also called sentential calculus or sentential logic)
is a formal system in which formulas of a formal language may be interpreted to represent propositions.

A system of rules and logical statements allows certain formulas to
be derived. These derived formulas may be interpreted to be true
propositions. Usually in Truth-functional propositional logic,
formulas are interpreted as having either a truth value of
\emph{true} or a truth value of \emph{false}.

Using the sentential logic in this paper we propose a composing
technology of new problems as an interpretation of  specific logical
models. Our aim is to give suitable logical models for formulation
of \emph{equivalent} problems and \emph{generating} problems of a
given problem.

\section{Preliminaries}

In logic, a set of symbols is commonly used to express logical representations.
Let us recall the basic symbols and logical representations we shall
deal with.

Let $p$ and $q$ be two statements.

\begin{itemize}
\item[i)] $p\,\wedge\,q \;$ denotes \emph{logical conjunction} (should be read as ``$p$ \,\emph{and}\, $q$'').
The statement $p\wedge q$ is \emph{true} if $p$ and $q$ are both true; else it is \emph{false}.
\vskip 2mm

\item[ii)] $p\, \vee\, q \;$ denotes \emph{logical disjunction} (should be read as ``$p$\, \emph{or} \,$q$'').
The statement $p\vee q$ is \emph{true} if $p$ or $q$ (or both) are true; if both are false,
the statement is \emph{false}.
\vskip 2mm

\item[iii)] $p\,\veebar\, q\;$ denotes \emph{exclusive disjunction} (should be read as ``\emph{either $p$ or $q$}'').
The statement $p\veebar q$ is \emph{true} when either $p$ or $q$, but not both, are true.
\vskip 2mm

\item[iv)] $\neg \,p\;$ denotes \emph{negation} (should be read as ``\emph{not}\, $p$'').
The statement $\neg p$ is \emph{true} if and only if $p$ is false.
\vskip 2mm

\item[v)] $p\, \rightarrow\, q\;$ denotes \emph{logical implication} (should be read as ``\emph{if $p$ then $q$}'').
The statement
$p \rightarrow q$ is \emph{true} just in the case that either $p$ is false or $q$ is true, or both.
The statements $p$ and $q$ aren't necessarily related comprehensively to each other.
\vskip 2mm

\item[vi)] $p\, \Rightarrow\, q\;$ denotes \emph{material implication} (should be read as ``\emph{$p$ implies $q$}'' or
``\emph{$q$ follows $p$}''). The relation
$p \Rightarrow q$ means that if $p$ is true then $q$ is also true; if $p$ is false then nothing is said about $q$.
The statements $p$ and $q$ are related comprehensively to each other.
\vskip 2mm

\item[vii)] $p\, \leftrightarrow\, q\;$ denotes \emph{logical equivalence} (should be read as
``$p\,$ \emph{if and only if}$\,q$''). The statement
$p \leftrightarrow q$ is \emph{true} just in case either both $p$ and $q$ are false, or both $p$ and $q$ are true.
The statements $p$ and $q$ aren't necessarily related comprehensively to each other.
\vskip 2mm

\item[viii)] $p\, \Leftrightarrow\, q\;$ denotes \emph{material equivalence}
(should be read as ``\emph{$q$ is necessary and sufficient for $p$''}). The relation
$p\,\Leftrightarrow\,q$ means that $p\,\Rightarrow\, q$ and $q\,\Rightarrow\, p$.
The statements $p$ and $q$ are related comprehensively to each other.
\end{itemize}

\section{Theoretical basis of the proposed method for generating problems}

In this section we describe in detail the theoretical basis of the method
for generating problems with a given logical structure.
In what follows $\;p_1,\, p_2;\, t,\, p,\, q,\, r\;$ will stand for statements.

In this paper we deal with a generalization of the formal logical rule \cite{S}
$$(p_1\rightarrow r)\,\wedge\, (p_2\rightarrow r)\;\Leftrightarrow \;
(p_1\vee p_2 \rightarrow r).\leqno(*)$$

Semantic rules connected with the material implication correspond to the formal derivation rules
used in the proofs below. By
semantic interpretations the formal derivation rules are called \emph{consequence rules} \cite{W}.

This correspondence allows us to formulate and comprehensively use the proposition below.
\vskip 2mm

\begin{prop}
The following equivalence is true
$$(t\wedge p\rightarrow r)\wedge (t\wedge q\rightarrow r)\;\Leftrightarrow\; t\wedge (p\vee q) \rightarrow r.
\leqno(1)$$
\end{prop}
\vskip 2mm

\emph{Proof}.
Let statement $p_1$ in (*) have structure $t\wedge p$ and  statement $p_2$ in (*) have  structure
$t\wedge q$. Then
$$(t\wedge p\rightarrow r)\wedge (t\wedge q\rightarrow r)\; \Leftrightarrow \;(t\wedge p)\vee (t\wedge q)\rightarrow r
\;\Leftrightarrow \;t\wedge (p\vee q) \rightarrow r,$$
i. e. the conjunction of the problems
$$t\wedge p\rightarrow r \leqno(2)$$
and
$$t\wedge q\rightarrow r \leqno(3)$$
is equivalent to the problem
$$t\wedge (p\vee q) \rightarrow r. \leqno(4)$$
\hfill{$\square$}
\vskip 2mm

Any true proposition could have more than one inverse proposition.
However, not every inverse proposition is
a true statement. The truth value of an inverse proposition of a given true proposition depends essentially
on its composition principle.

According to \cite{S}, if a given proposition has the logical structure
$\;p_1\,\wedge\, p_2\,\rightarrow\, r$, then each one of the following propositions could be considered to be its
inverse:\\ $\;r\,\rightarrow\,p_1\,\wedge\, p_2$, $\;p_1\,\wedge \,r\,\rightarrow\, p_2\;$ and
$\;r\,\wedge \,p_2\,\rightarrow\, p_1$.

The most interesting and important inverse propositions are those
that are true as well as independent from the other possible inverse propositions, i. e. the \emph{strongest}
inverse propositions.

Equivalence (1) formally describes a method for composing new problems with a given logical structure
and for formulating their inverse problems.

According to Proposition 3.1 problems with logical structures (2) and (3) generate a problem with
a logical structure (4).

In this paper we consider only problems inverse to problems of type (4) with the structure
$$t\,\wedge\,r\;\rightarrow\;p\,\vee\,q. \leqno(5)$$

Problems with logical structures (2) and (3) are said to be \emph{generating} problems with
structure (4) and their inverse problems with structure (5).

To change the logical structure in the conclusion of the inverse
problem from \emph{logical disjunction} to \emph{exclusive
disjunction} we need a dichotomic decomposition of the considered
set of geometric objects with respect to any remarkable property and
its negation. Such a decomposition guarantees the \emph{homogeneity} \cite{P}
of the statements (based on one and the same equivalence relation)
in the conclusion of the problem.
\vskip 2mm

\begin{prop}
The following equivalence is true
$$t\,\wedge\,r\;\rightarrow\;p \,\vee\, q\;\Leftrightarrow\;
t\,\wedge\,r\;\rightarrow\;p \,\veebar\,(\neg p\,\wedge\,q). \leqno(6)$$
\end{prop}
\vskip 2mm

Proposition 3.2 gives the equivalence between problems with a logical structure (5) and problems
with a logical structure
$$t\,\wedge\,r\;\rightarrow\;p \,\veebar\,(\neg p\,\wedge\,q). \leqno(7)$$

Any problem with a logical structure (7) satisfies the condition of \emph{homogeneity in the conclusion}.

\section{Application of the method to specific groups of problems}

We discuss four groups of problems to illustrate the described generating method.
In each of the groups we formulate  suitable \emph{generating} problems for
the corresponding equivalent and inverse problems.

The problems in each of the proposed groups are comprehensively related to each other.

\subsection{Problems of group I}

The statements used for the formulation of the problems in this
group are

\begin{itemize}
\item $t$:=\{\it The straight line $AD,\, D\in BC$, is a median in  $\triangle \,ABC$.\}
\vskip 2mm

\item $p$:= \{$AC=AB$\}
\vskip 2mm

\item $q$:=\{$\angle BAC=90^0$\}
\vskip 2mm

\item $r$:=\{$\angle DAC + \angle ABC = 90^0$\}
\end{itemize}
\vskip 2mm

First we formulate  the \emph{generating} problems.
\begin{prob}
Let the straight line $AD,\, D\in BC$, be a median in  $\triangle\, ABC$. Prove that if  $AC=AB$,
then $\angle DAC + \angle ABC = 90^0$.
\end{prob}

\begin{figure}[h t b]
\epsfxsize=8cm
\centerline{\epsfbox{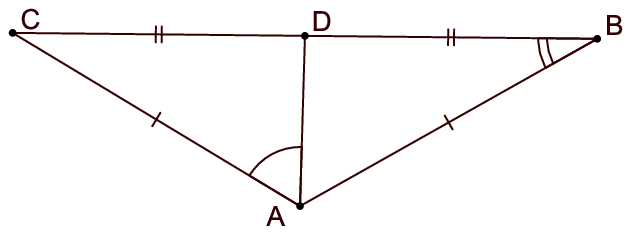}}
\end{figure}

This problem has a logical structure $t\,\wedge\, p \;\rightarrow\; r$.
Its proof follows immediately from Figure 1.
\vskip 2mm

\begin{prob}
Let the straight line $AD,\, D\in BC$, be a median in  $\triangle\, ABC$. Prove that if  $\angle BAC=90^0$,
then $\angle DAC + \angle ABC = 90^0$.
\end{prob}

Problem 4.2 has a logical structure $t\,\wedge\, q \;\rightarrow\; r$.
The proof follows easily from Figure 2.

\begin{figure}[h t b]
\epsfxsize=8cm
\centerline{\epsfbox{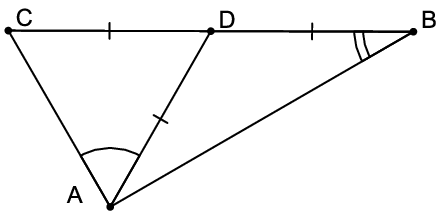}}
\end{figure}
\vskip 2mm

According to the logical structures of Problems 4.1 and 4.2 and in view of Proposition 3.1, we construct
the following inverse problem with logical structure
$\,t\,\wedge\,r\;\rightarrow\;p\,\vee\,q\,$.
\vskip 2mm

\begin{prob} $($\cite{NM}, Problem $3)$
Let the straight line $AD,\, D\in BC$, be a median in  $\triangle \,ABC$.
Prove that if $\,\angle DAC + \angle ABC = 90^0$, then
$\,AC=AB\,$ or $\,\angle BAC=90^0$.
\end{prob}
\vskip 2mm

The next two problems are equivalent to Problem 4.3.
\vskip 2mm

\begin{prob} $($\cite{NM}, Problem $2)$
Let the straight line $AD,\, D\in BC$, be a median in  $\triangle \,ABC$.
Prove that if $\,\angle DAC + \angle ABC = 90^0$
and $\angle BAC\neq 90^0$, then $\,AC=AB$.
\end{prob}
\vskip 2mm

\begin{prob} $($\cite{NM}, Problem $1)$
Let the straight line $AD,\, D\in BC$, be a median in  $\triangle \,ABC$.
Prove that if $\,\angle DAC + \angle ABC = 90^0$
and $\,AC\neq AB$, then $\angle BAC = 90^0$.
\end{prob}

Another version of Problem 4.5 is Problem 246, p. 211 in \cite{Sch3}.
\vskip 2mm

In view of Proposition 3.2, Problem 4.3 can be reformulated as follows by keeping the condition of
homogeneity in its conclusion (compare also with \cite{NM},Problem 4;
\cite{Sch2}, p. 24, Problem 6; \cite{Sch1}, p. 22, Problem 1; \cite{Sch4}, p. 265, Problem 312).
\vskip 2mm

\begin{prob}
Let the straight line $AD,\, D\in BC$, be a median in  $\triangle \,ABC$.
Prove that if $\,\angle DAC + \angle ABC = 90^0$,
then $\triangle \,ABC$ is either isosceles $\,(AC = AB)$, or not isosceles but right-angled $(\angle BAC = 90^0)$.
\end{prob}
\vskip 4mm

\subsection{Problems of group II}

The statements used for the formulation of the problems in this
group are

\begin{itemize}
\item $t$:=\{\it In $\triangle \,ABC$ the straight line $AA_1,\,A_1\in BC$, is the bisector of
$\angle CAB$, the straight line
$BB_1,\, B_1\in AC$, is the bisector of $\angle CBA$ and $AA_1\cap BB_1=J$.\}
\vskip 2mm

\item $p:= \{AC=BC\}$
\vskip 2mm

\item $q:=\{\angle ACB=60^0\}$
\vskip 2mm

\item $r:=\{JA_1=JB_1\}$
\end{itemize}
\vskip 2mm

First we formulate and solve the \emph{generating} problems.

\begin{prob}
Let in $\triangle \,ABC$ the straight line $AA_1,\,A_1\in BC$, be the bisector of
$\angle CAB$, the straight line
$BB_1,\, B_1\in AC$, be the bisector of $\angle CBA$ and $AA_1\cap BB_1=J$.
Prove that if $AC=BC$, then $JA_1=JB_1$.
\end{prob}

Problem 4.7 has a logical structure $t\,\wedge\, p \;\rightarrow\; r$.

\begin{figure}[h t b]
\epsfxsize=8cm
\centerline{\epsfbox{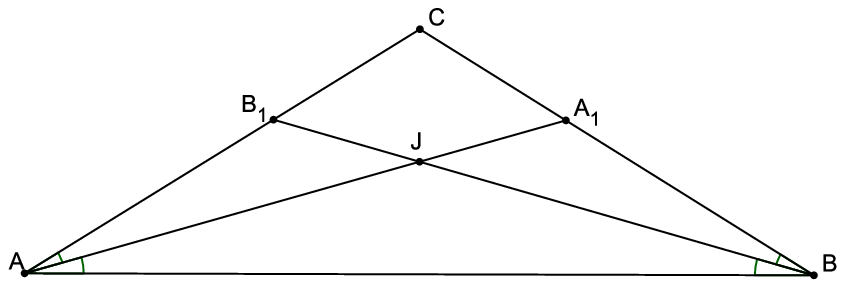}}
\end{figure}

\emph{Proof}. Since $AC=BC$, then $\angle CAB=\angle CBA$ and hence $\angle A_1AB=\angle B_1BA$ (fig. 3).

From the Criteria  for congruence of triangles we have $\,\triangle A_1AB\cong\triangle B_1BA$. As a consequence
it follows that  $\,AA_1=BB_1,$  $\,\triangle AJB\,$ is isosceles, $AJ=BJ\;$ and $\, JA_1=JB_1$.

\hfill{$\square$}
\vskip 2mm

\begin{prob}
Let in $\triangle \,ABC$ the straight line $AA_1,\,A_1\in BC$, be the bisector of
$\angle CAB$, the straight line
$BB_1,\, B_1\in AC$, be the bisector of $\angle CBA$ and $AA_1\cap BB_1=J$.
Prove that if $\angle ACB=60^0$, then $JA_1=JB_1$.
\end{prob}

This problem has a logical structure $t\,\wedge\, q \;\rightarrow\; r$.
\vskip 2mm

\begin{figure}[h t b]
\epsfxsize=8cm
\centerline{\epsfbox{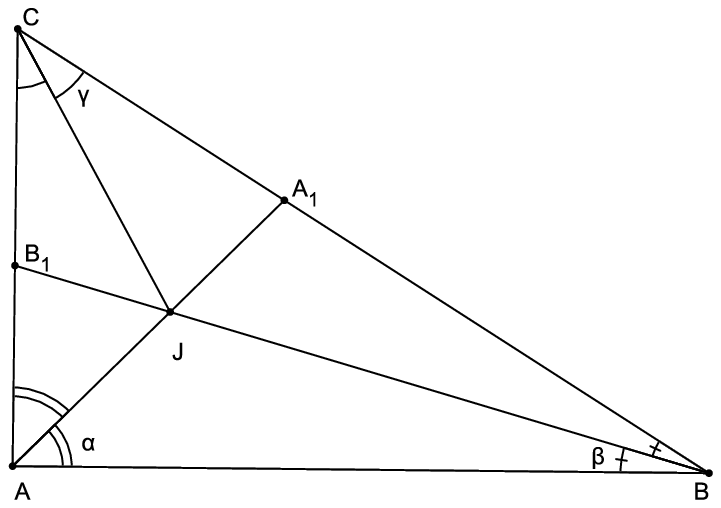}}
\end{figure}

\emph{Proof}. Let us denote $\,\angle BAA_1=\angle CAA_1=\alpha,\, \angle ABB_1=\angle CBB_1=\beta$ (fig. 4).

Since $J$ is the intersection point of the bisectors $AA_1$ and $BB_1$ of $\triangle \,ABC$, then $CJ$ is the
bisector of $\angle ACB$ and $\angle JCA=\angle JCB=\gamma=30^0$. Since $\alpha+\beta+\gamma=90^0$, then
$\alpha+\beta=60^0$, $\angle AJB=120^0$ and the quadrilateral $CA_1JB_1$ can be inscribed in a circle. Hence,
$JA_1=JB_1$ as chords corresponding to equal angles (arcs) in a circle.
\hfill{$\square$}
\vskip 2mm

According to the logical structures of Problems 4.7 and 4.8 and in view of Proposition 3.1 we construct
the following inverse problem with logical structure
$\,t\,\wedge\,r\;\rightarrow\;p\,\vee\,q\,$.
\vskip 2mm

\begin{prob} $($\cite{NM}, Problem $7)$
Let in $\triangle \,ABC$ the straight line $AA_1,\,A_1\in BC$, be the bisector of
$\angle CAB$, the straight line
$BB_1,\, B_1\in AC$, be the bisector of $\angle CBA$ and $AA_1\cap BB_1=J$.
Prove that if $JA_1=JB_1$, then
$AC=BC$ or $\angle ACB=60^0$.
\end{prob}
\vskip 2mm

The next two problems are equivalent to Problem 4.9.
\vskip 2mm

\begin{prob} $($\cite{NM}, Problem $5)$
Let in $\triangle \,ABC$ the straight line $AA_1,\,A_1\in BC$, be the bisector of
$\angle CAB$, the straight line
$BB_1,\, B_1\in AC$, be the bisector of $\angle CBA$ and $AA_1\cap BB_1=J$.
Prove that if $JA_1=JB_1$ and
$AC\neq BC$,  then $\angle ACB=60^0$.
\end{prob}
\vskip 2mm

\begin{prob} $($\cite{NM}, Problem $6)$
Let in $\triangle \,ABC$ the straight line $AA_1,\,A_1\in BC$, be the bisector of
$\angle CAB$, the straight line
$BB_1,\, B_1\in AC$, be the bisector of $\angle CBA$ and $AA_1\cap BB_1=J$.
Prove that if $JA_1=JB_1$ and
$\angle ACB \neq60^0$, then $AC=BC$.
\end{prob}
\vskip 2mm

In view of Proposition 3.2, Problem 4.9 can be reformulated by keeping the condition of
homogeneity in its conclusion.
\vskip 2mm

\begin{prob} $($\cite{NM}, Problem $8)$
Let in $\triangle \,ABC$ the straight line $AA_1,\,A_1\in BC$, be the bisector of
$\angle CAB$, the straight line
$BB_1,\, B_1\in AC$, be the bisector of $\angle CBA$ and $AA_1\cap BB_1=J$.
Prove that if $JA_1=JB_1$
then $\triangle \,ABC$ is either isosceles $\,(CA = CB)$, or not isosceles but $\angle ACB = 60^0$.
\end{prob}
\vskip 4mm

\subsection{Problems of group III}

Let in $\triangle \,ABC$ the straight line $CH,\,H\in AB$, be the
altitude and the straight line $CM,\, M\in AB$, be the median.

In  the not isosceles $\,\triangle \,ABC\,$ the location of the
collinear points $A$, $B$, $H$ and $M$ is either $M/AH$ or $M/BH$.
If the triangle is isosceles, then $M\equiv H$.

Denoting as usually $\,\angle CAB\,$ by $\,\alpha$ and $\,\angle
CBA\,$ by $\,\beta$, exactly one of the following relations is
fulfilled: $\,\alpha < \beta\;$ (in this case $\,M/AH$), $\quad
\alpha = \beta\;$ (in this case $\,M\equiv H$), $\quad\alpha >
\beta\;$ (in this case $M/BH$).

In the case $\alpha < \beta \;(M/AH)$ we compare the angles
$\angle\, ACH$ and $\angle\, BCM$; in the case $\alpha > \beta
\;(M/BH)$ we compare the angles $\angle\, ACM$ and $\angle\, BCH$.
In both cases the considerations are analogical to each other.
What is more, if both angles $\alpha$ and $\beta$ are acute
angles, i. e. $H/AB$, then
$$\angle\, ACH=\angle\, BCM\quad\Leftrightarrow\quad\angle\, ACM=\angle\, BCH.$$
\vskip 2mm

In what follows let $\alpha \geq \beta$, i. e. either $\,M/BH$, or
$\,M\equiv H$. \vskip 1mm

The statements used for the formulation of the problems in this
group are

\begin{itemize}
\item $t:=\{${\it The straight line $CH,\,H\in AB$, is the altitude and the straight
line $CM,\, M\in AB$, is the median  of $\triangle \,ABC$.}$\}$
\vskip 2mm

\item $p:=\{ AC=BC\}$
\vskip 2mm

\item $q:=\{\angle ACB= 90^0\}$
\vskip 2mm

\item $r:=\{\angle ACM=\angle BCH\}\;$
\end{itemize}
\vskip 2mm

First we formulate and solve two \emph{generating} problems.

\begin{prob}
Let the straight line $CH,\,H\in AB$, be the altitude and the
straight line $CM,\, M\in AB$, be the median  of $\triangle
\,ABC$. If $AC=BC$, then $\angle ACM=\angle BCH$.
\end{prob}

This problem has a logical structure $t\,\wedge\, p
\;\rightarrow\; r$. \vskip 2mm

\emph{Proof}. In any isosceles triangle the altitude and the
median to its base are congruent. Hence, $M\equiv H$ and $\angle
ACM=\angle BCH$. \hfill{$\square$} \vskip 2mm

\begin{prob}
Let the straight line $CH,\,H\in AB$, be the altitude and the
straight line $CM,\, M\in AB$, be the median  of $\triangle
\,ABC$. If $\angle ACB= 90^0$, then $\angle ACM=\angle BCH$ $($and
also $\angle ACH=\angle BCM)$.
\end{prob}

This problem has a logical structure $t\,\wedge\, q
\;\rightarrow\; r$. \vskip 2mm

\emph{Proof}. In the right-angled not isosceles $\triangle \,ABC$
the location of the collinear points $B$, $H$ and $M$ is either
$H/BM$, or $H/AM$. In the case under consideration $\alpha >\beta$
and $H/AM$ (fig. 5).
\begin{figure}[h t b]
\epsfxsize=8cm \centerline{\epsfbox{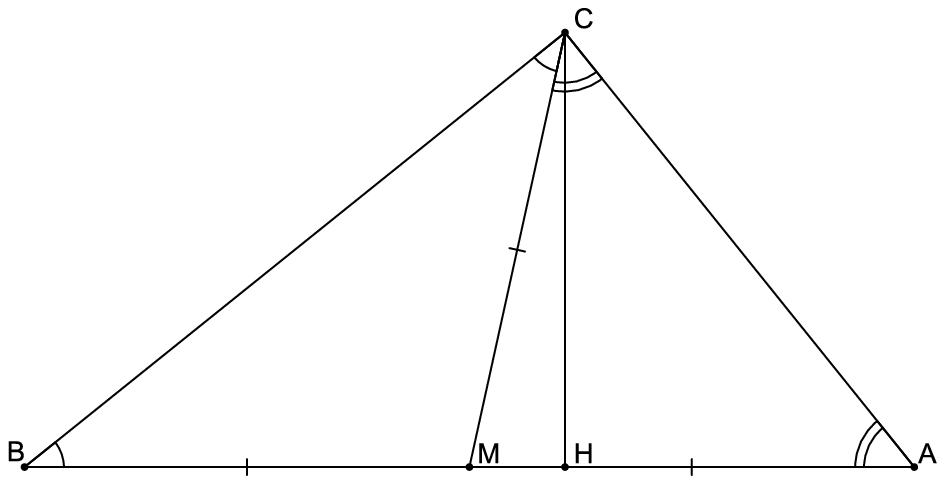}}
\end{figure}
On the other hand  $AM=MC\;(=MB)$ and $\triangle AMC$ is
isosceles. Thus, $\angle ACM=\alpha$. Since $\alpha+\beta=90^0$,
in the right-angled $\triangle BHC$ we have  $\angle
BCH=90^0-\beta=\alpha$. Hence, $\angle ACM=\angle BCH$ (and also
$\angle ACH=\angle BCM$).

For a right-angled  isosceles triangle see Problem 4.13.
\hfill{$\square $}

\vskip 2mm According to the logical structures of problems 4.13
and 4.14 and in view of Proposition 3.1 we construct the following
inverse problem with logical structure
$\,t\,\wedge\,r\;\rightarrow\;p\,\vee\,q\,$. \vskip 2mm

\begin{prob}
Let the straight line $CH,\,H\in AB$, be the altitude and the
straight line $CM,\, M\in AB$, be the median  of $\triangle
\,ABC$. If $\angle ACM=\angle BCH$, then $AC=BC$ $($i. e.
$\triangle \,ABC$ is isosceles$)$ or $\angle ACB= 90^0$ $($i. e.
$\triangle \,ABC$ is right-angled$)$.
\end{prob}
\vskip 2mm

\emph{Proof}.
\begin{figure}[h t b]
\epsfxsize=8cm \centerline{\epsfbox{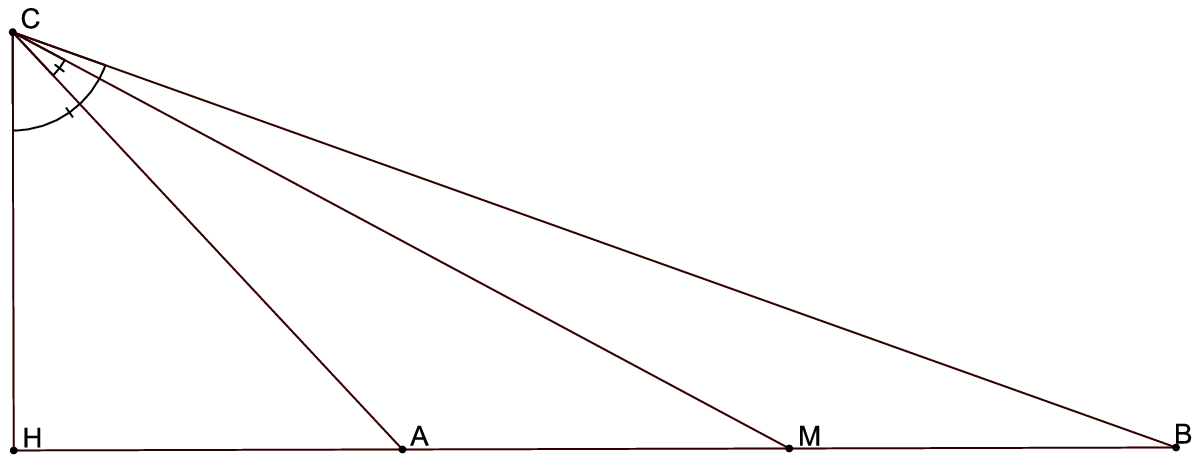}}
\end{figure}
In any triangle at least two of the angles must be acute angles.
Hence, in $\triangle \,ABC$ at least one of the angles $\alpha$
and $\beta$ is acute. Since $\alpha \geq \beta$, then  $\beta<
90^0$.

If we assume that $\alpha\geq 90^0$ then the location of the
collinear points $A$, $H$ and $M$ is either $A/HM$, or $A\equiv H$
(fig. 6). Then for the right-angled $\triangle BCH$ is valid
$\angle ACM<\angle BCH$, which contradicts the given condition
$\angle ACM=\angle BCH$. Hence, $\alpha<90^0$ and the points $H$
and $M$ lie between the points $A$ and $B$.

There are two possibilities for the points $H$ and $M$ -  they
either coincide or not. \vskip 2mm

$(i)$ Let $H\equiv M$. In this case the median $CM$ in $\triangle
\,ABC$ coincides with the altitude $CH$, i. e. $\triangle \,ABC$
is isosceles. If in addition $\angle ACB=90^0$, then $\triangle
\,ABC$ is isosceles  right-angled. \vskip 2mm

$(ii)$ Let $H\neq M$. Since $\alpha\geq \beta$, then $H/AM$ (fig.
7).
\begin{figure}[h t b]
\epsfxsize=8cm \centerline{\epsfbox{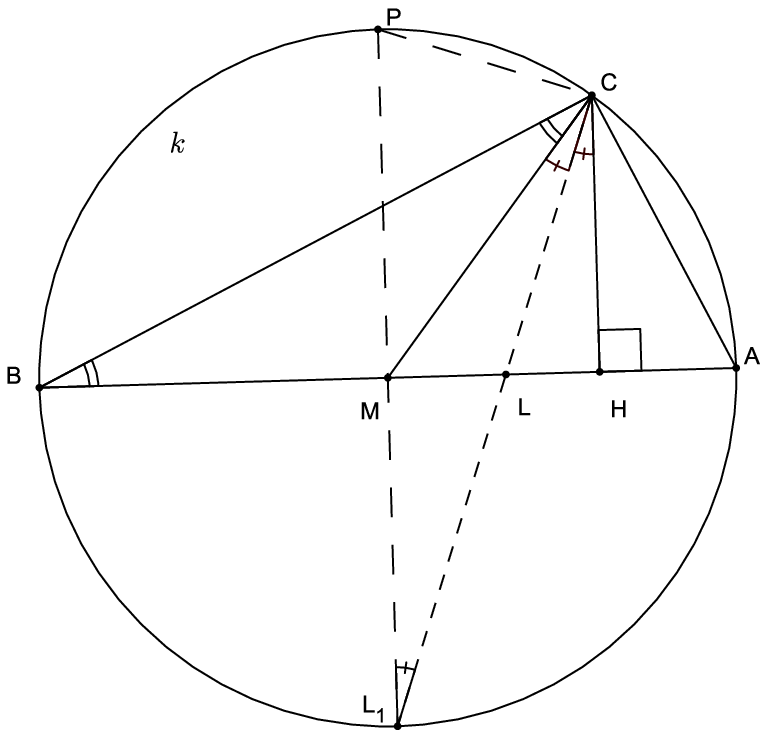}}
\end{figure}

Let $CL,\, L\in AB$, be the bisector of $\angle ACB$. It follows
that $CL$ is also the bisector of $\angle MCH$ (see also \cite{M},
p. 184, problem 29; \cite{Pr}, p. 41, problem 2.32).

Let $k$ be the circumscribing circle of $\triangle \,ABC$ and
$CL\cap k=L_1$. The point $L_1$ is the middle point of the arc
$\widehat{AL_1B}$. The points $C$ and $L_1$ lie on alternate sides
of $AB$. The perpendicular projection of $L_1$ onto the chord $AB$
is the middle point $M$. Then the straight line $L_1M$ is the
perpendicular bisector of $AB$.

The straight line $CL_1$ cuts the parallel lines $CH\, (CH\perp
AB)$ and $L_1M\, (L_1M\perp AB)$ and hence the alternate angles
$\angle HCL$ and $\angle ML_1L$ are equal, i. e. $\triangle CML_1$
is isosceles. Thus the point $M$ also lies on the the
perpendicular bisector of the chord $CL_1$.

Since the perpendicular bisectors of any two non parallel chords
of a circle cut at its center, the point $M$ is the center of $k$,
the chord $AB$ is a diameter of $k$ and $\angle ACB=90^0$.

\begin{rem}
Let $P=ML_1\cap k$. Then $PL_1$ is a diameter of $k$ and $\angle
PCL_1=90^0$. It is easy to be seen that $\triangle MPC$ is
isosceles and the point $M$ is the center  of $k$.
\end{rem}
\hfill{$\square$} \vskip 2mm

We reformulate Problem 4.15 by keeping the condition of
homogeneity of the conclusion. \vskip 2mm

\begin{prob}
Let the straight line $CH,\,H\in AB$, be the altitude and the
straight line $CM,\, M\in AB$, be the median  of $\triangle
\,ABC$. If $\angle ACM=\angle BCH$, then $\triangle ABC$ is either
isosceles $(AC=BC)$, or not isosceles but right-angled $(\angle
ACB=90^0)$.
\end{prob}
\vskip 4mm

\subsection{Problems of group IV}

The statements used for the formulation of the problems in this
group are

\begin{itemize}
\item $t:=\{$\it  The middle points of the sides $BC$, $CA$ and $AB$ of $\triangle ABC$ are
$F$, $D$, and $E$ respectively.$\}$

\vskip 2mm
\item $p:=\{AC=BC\}$
\vskip 2mm

\item $q:=\{\angle ACB=60^0\}$
\vskip 2mm

\item $r:=\{$ The center $G$ of the circumscribing circle $k$ of $\triangle FDE$ lies on the bisector
of $\angle ACB\}.$
\end{itemize}
\vskip 2mm

First we formulate and solve the \emph{generating} problems.

\begin{prob}
Let the middle points of the sides $BC$, $CA$ and $AB$ of $\triangle ABC$ be
$F$, $D$, and $E$ respectively. Prove that if $AC=BC$, then the center $G$ of the circumscribing
circle $k$ of $\triangle FDE$ lies on the bisector of $\angle ACB$.
\end{prob}

This problem has a logical structure $t\,\wedge\, p \;\rightarrow\; r$.
\vskip 2mm

\emph{Proof}. The median $CE$ of the isosceles $\triangle ABC$ is the perpendicular bisector of
$AB$ and $DF$ and the bisector of $\angle ACB$. Hence, the center $G$ of the circumscribing
circle $k$ of $\triangle FDE$ lies on the bisector of $\angle ACB$.
\hfill{$\square$}
\vskip 2mm

\begin{prob}
Let the middle points of the sides $BC$, $CA$ and $AB$ of $\triangle ABC$ be
$F$, $D$, and $E$ respectively.
Prove that if $\angle ACB=60^0$, then the center $G$ of the circumscribing
circle $k$ of $\triangle FDE$ lies on the bisector of $\angle ACB$.
\end{prob}

\begin{figure}[h t b]
\epsfxsize=8cm
\centerline{\epsfbox{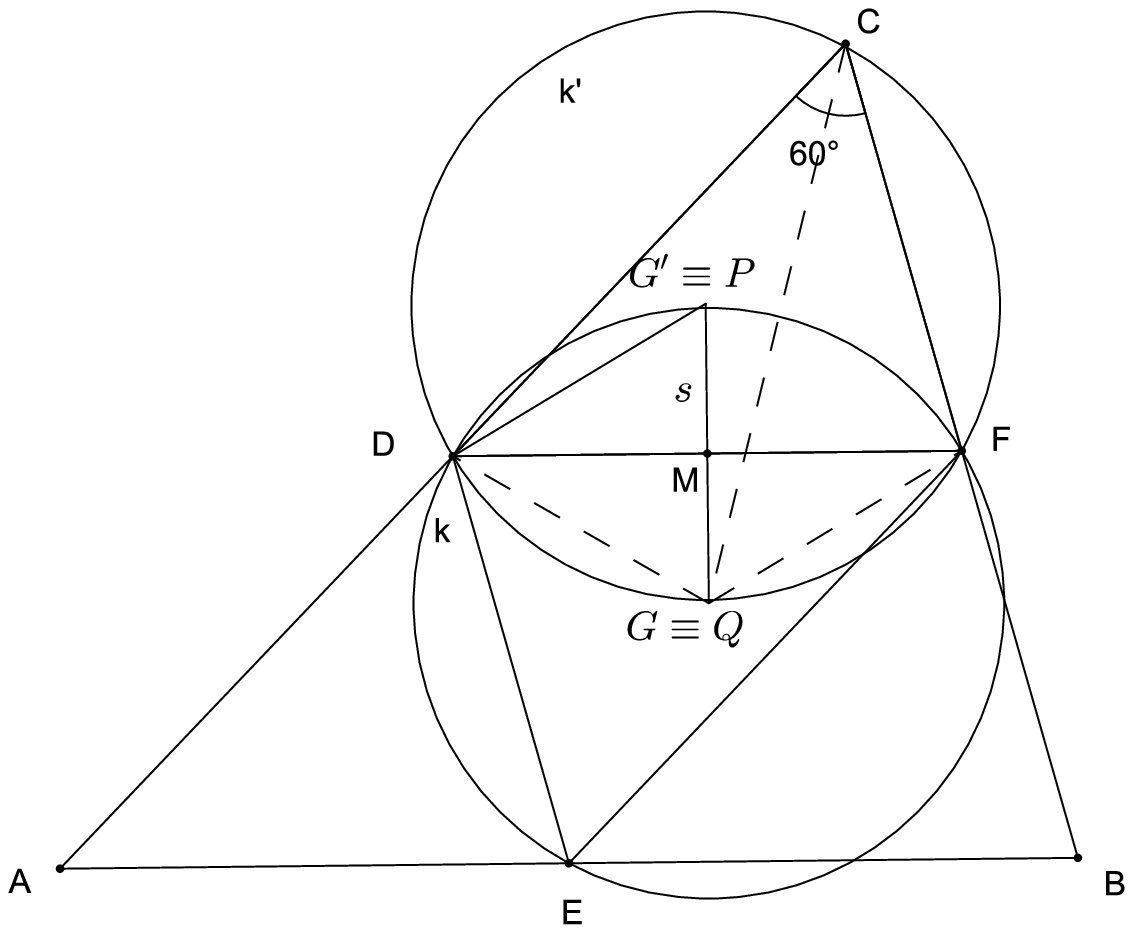}}
\end{figure}

This problem has a logical structure $t\,\wedge\, q \;\rightarrow\; r$.
\vskip 2mm

\emph{Proof}.
The quadrilateral $EFCD$ (fig. 8) is a parallelogram with $\angle DCF=60^0$.
Hence, $\triangle EFD\cong\triangle CDF$ and the circumscribing circles $k$  and $k'$ of
$\triangle EFD$ and $\triangle CDF$ respectively have equal radii. The centers $G$ and $G'$
of these circles lie on the perpendicular bisector $s$ of $DF$.

Let $P=s\cap k,\;$ $Q=s\cap k'$.
It is easy to be seen that the quadrilateral $FPDQ$ is a rhombus with $\angle PDQ=60^0$
and $QD=QP=QF$, i. e. the point $Q$ coincides with the center $G$ of $k$. Consequently, the
point $P$ coincides with the center $G'$ of $k'$.

The point $Q$
is also the middle point of the arc $\widehat{DQF}$ of $k'$ and then lies on the bisector
of $\angle DCF\equiv \angle ACB$.
\hfill{$\square$}
\vskip 2mm

According to the logical structures of Problems 4.18 and 4.19 and in view of Proposition 3.1 we construct
the following inverse problem with logical structure $\,t\,\wedge\,r\;\rightarrow\;p\,\vee\,q\,$
(a formulation with a different logical structure is given in  \cite{T}, Problem 12):
\vskip 2mm

\begin{prob}
Let the middle points of the sides $BC$, $CA$ and $AB$ of $\triangle ABC$ be
$F$, $D$, and $E$ respectively.
Prove that if the center $G$ of the circumscribing
circle $k$ of $\triangle FDE$ lies on the bisector of $\angle ACB$, then
$AC=BC$ or $\angle ACB=60^0$.
\end{prob}
\vskip 2mm

The next two problems are equivalent to Problem 4.20.
\vskip 2mm

\begin{prob}
Let the middle points of the sides $BC$, $CA$ and $AB$ of $\triangle ABC$ be
$F$, $D$, and $E$ respectively.
Prove that if the center $G$ of the circumscribing
circle $k$ of $\triangle FDE$ lies on the bisector of $\angle ACB$ and $BC\neq AC$, then
$\angle ACB=60^0$.
\end{prob}
\vskip 2mm

\begin{prob}
Let the middle points of the sides $BC$, $CA$ and $AB$ of $\triangle ABC$ be
$F$, $D$, and $E$ respectively.
Prove that if the center $G$ of the circumscribing
circle $k$ of $\triangle FDE$ lies on the bisector of $\angle ACB$ and $\angle ACB\neq 60^0$,
then $BC=AC$.
\end{prob}
\vskip 2mm

We reformulate Problem 4.20 by keeping the condition of homogeneity of the conclusion.
\vskip 2mm

\begin{prob}
Let the middle points of the sides $BC$, $CA$ and $AB$ of $\triangle ABC$ be
$F$, $D$, and $E$ respectively.
Prove that if the center $G$ of the circumscribing
circle $k$ of $\triangle FDE$ lies on the bisector of $\angle ACB$, then
the $\triangle ABC$ is either isosceles $(AC=BC)$, or not isosceles but $\angle ACB=60^0$.
\end{prob}

\begin{figure}[h t b]
\epsfxsize=8cm
\centerline{\epsfbox{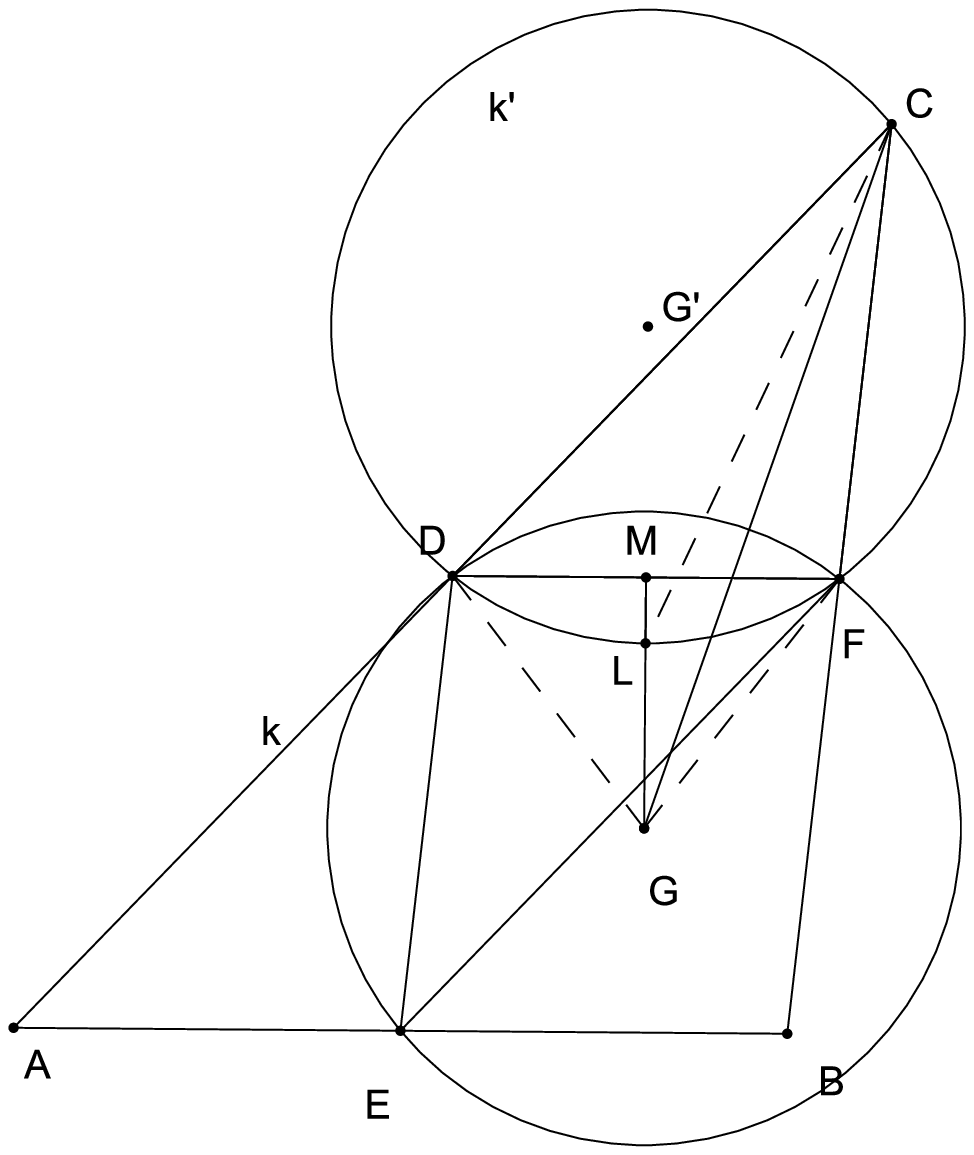}}
\end{figure}

\emph{Proof}.  Let $G'$ be the center of the circumscribing circle $k'$ of $\triangle FDC$ (fig. 9).
In view of the Criteria for congruence of triangles we get that $\triangle FDE\cong \triangle DFC$. It follows that
the circumscribing circles $k$  and $k'$ of
$\triangle FDE$ and $\triangle DFC$ respectively have equal radii.

Let $M$ be the middle point of
$DF$ and $L=GM\cap k'$. The point $G'$ lies on the perpendicular bisector $GM$ of $DF$.
Hence, the point $L$ is the middle point of the
arc $\widehat{DLF}$ of $k'$ and $CL$ is the bisector of $\angle DCF\equiv \angle ACB$.

Since the center $G$ of $k$ lies on the bisector $CL$ (according to the condition of the Problem),
then the straight lines $CL$ and $GM$
either cut at $G$ (have no other common points), or coincide (all of their points are common).
\vskip 2mm

$(i)$ Let $CL\cap GM=L\equiv G$.

In this case $G'\in k$ (fig. 8) and $\triangle G'DG$ is equilateral,
the central $\angle DG'F$ of $k'$ has a measure $120^0$ and hence $\angle ACB=60^0$.
\vskip 2mm

$(ii)$ Let $CL\equiv GM$ (fig. 10).
\begin{figure}[h t b]
\epsfxsize=8cm
\centerline{\epsfbox{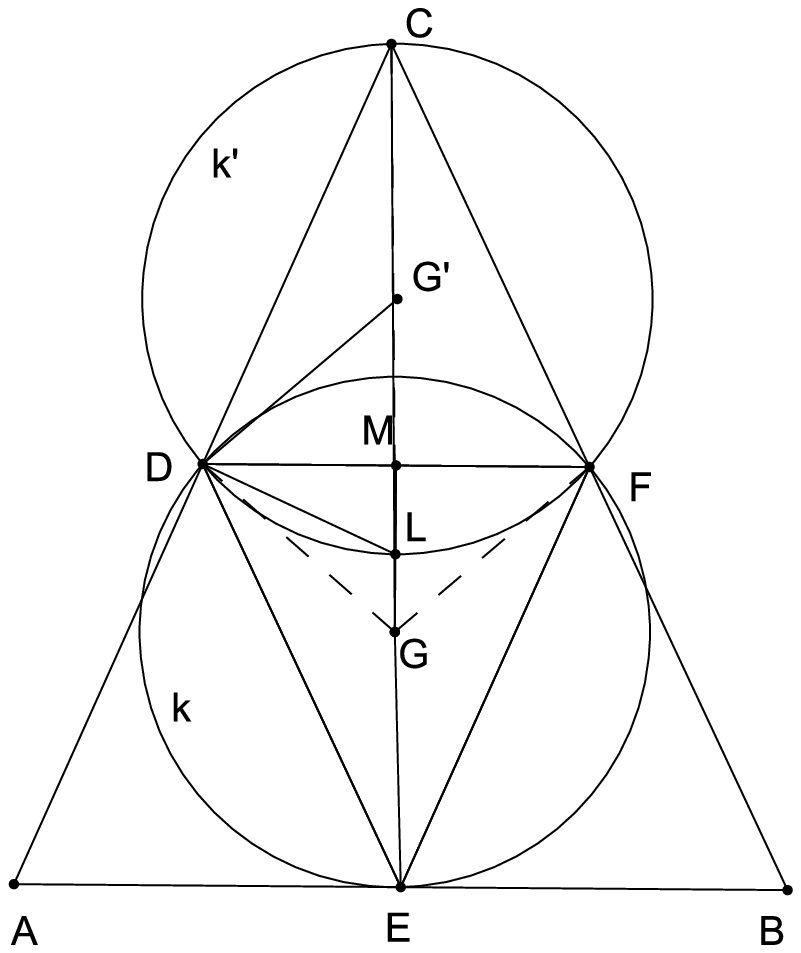}}
\end{figure}

In this case  the bisector $CL$ of $\angle DCF$ coincides with the perpendicular bisector of $DF$. Then
$\triangle DCF$ and also $\triangle ABC$ are isosceles, i. e. $AC=BC$.
\hfill{$\square$}
\vskip 6mm

\section{Summary}
In this section we formulate a new problem whose proof emphasizes the
importance and significance of the described method for
generating problems.
\vskip 2mm

The similar conclusions of Problems 4.9 and 4.20 lead to
\vskip 2mm

\begin{prob}
Let the middle points of the sides $\,BC$, $\,CA\,$ and $\,AB\,$ of $\,\triangle\, ABC\,$ be
$\,F$, $\,D\,$ and $\,E\,$ respectively. Let further the straight lines $\,AA_1,\,A_1\in BC$, and
$\,BB_1,\, B_1\in AC$, be the bisectors of $\,\angle \,CAB\,$ and $\,\angle CBA$, respectively,
and let $\,AA_1\cap BB_1=J$.

Prove that the center $\,G\,$ of the circumscribing
circle $\,k\,$ of $\,\triangle\, FDE\,$ lies on the bisector of $\,\angle\, ACB\,$
if and only if $\,JA_1=JB_1$.
\end{prob}
\vskip 2mm

\emph{Proof}.

(i) Let the center $G$ of the circumscribing
circle $k$ of $\,\triangle FDE\,$ lie on the bisector of $\angle ACB$.

From Problem 4.20 it follows that $AC=BC$ or $\angle ACB=60^0$.

\begin{itemize}
\item[-] If $AC=BC$ then from the generating Problem 4.7 it follows that $JA_1=JB_1$.
\vskip 2mm

\item[-] If $\angle ACB=60^0$ then from the generating Problem 4.8 it follows that $JA_1=JB_1$.
\end{itemize}
\vskip 2mm

(ii) Let $JA_1=JB_1$.
From Problem 4.9 it follows that $AC=BC$ or $\angle ACB=60^0$.

\begin{itemize}
\item[-] If $AC=BC$ then from the generating Problem 4.18 it follows that the center $G$ of the circumscribing
circle $k$ of $\triangle FDE$ lies on the bisector of $\angle ACB$.

\vskip 2mm

\item[-] If $\angle ACB=60^0$ then from the generating Problem 4.19 it follows that the center $G$ of the
circumscribing circle $k$ of $\triangle FDE$ lies on the bisector of $\angle ACB$.
\end{itemize}
\hfill{$\square$}
\vskip 6mm

Acknowledgements. The first author is partially supported by Sofia University Grant 159/2013.
The second author is partially supported by Sofia University Grant 99/2013.

\end{document}